\documentclass[11pt]{article}
\usepackage[english]{babel}
\usepackage[utf8]{inputenc}
\usepackage{amsmath,amsfonts,amssymb,amsthm,mathrsfs}
\usepackage{minitoc}
\usepackage{enumitem} 

\renewcommand{\geq}{\geqslant}
\usepackage[mono=false]{libertine}
\useosf
\usepackage{mathpazo}

\usepackage{wasysym}

\usepackage[dvipsnames]{xcolor}
\usepackage[a4paper,vmargin={3.5cm,3.5cm},hmargin={2cm,2cm}]{geometry}
\usepackage[font=sf, labelfont={sf,bf}, margin=1cm]{caption}
\usepackage{graphicx,graphics}
\usepackage{epsfig}
\usepackage{latexsym}
\usepackage{stmaryrd}
\usepackage{ae,aecompl}

\usepackage[english]{babel}
 \usepackage[colorlinks=true]{hyperref}
\usepackage{pstricks}
\usepackage{enumerate}

\DeclareFixedFont{\beaupetit}{T1}{ftp}{b}{n}{2cm}

\newtheorem{theorem}{Theorem}[]

\theoremstyle{definition}

\usepackage{todonotes}

\title{{\bf Last Car Decomposition of Planar Maps }}
\author{Alice Contat} 
             
\begin{document}
\date{ }            
             \maketitle
\abstract{We give new equations which characterize the generating functions of planar quadrangulations and planar triangulations, with zero, one or two boundaries. The proof is inspired by the Lackner--Panholzer last car decomposition of parking trees \cite{LaP16} and consists in applying a similar decomposition to the peeling trees of planar maps.}

\section{Introduction}
 
Since their introduction by Tutte in his series of ``census" papers \cite{Tut62c,Tut62,Tut62b,Tut63}, maps are fundamental objects which have been extensively studied especially in combinatorics and in probability. The purpose of this work is to establish new recursive equations to enumerate particular types of maps which are inspired both by the well-known peeling procedure for maps \cite{CurStFlour} and by the last car decomposition of fully parked trees introduced first by Lackner and Panholzer in \cite{LaP16}.

\begin{figure}[!h]
 \begin{center}
 \includegraphics[width=15cm]{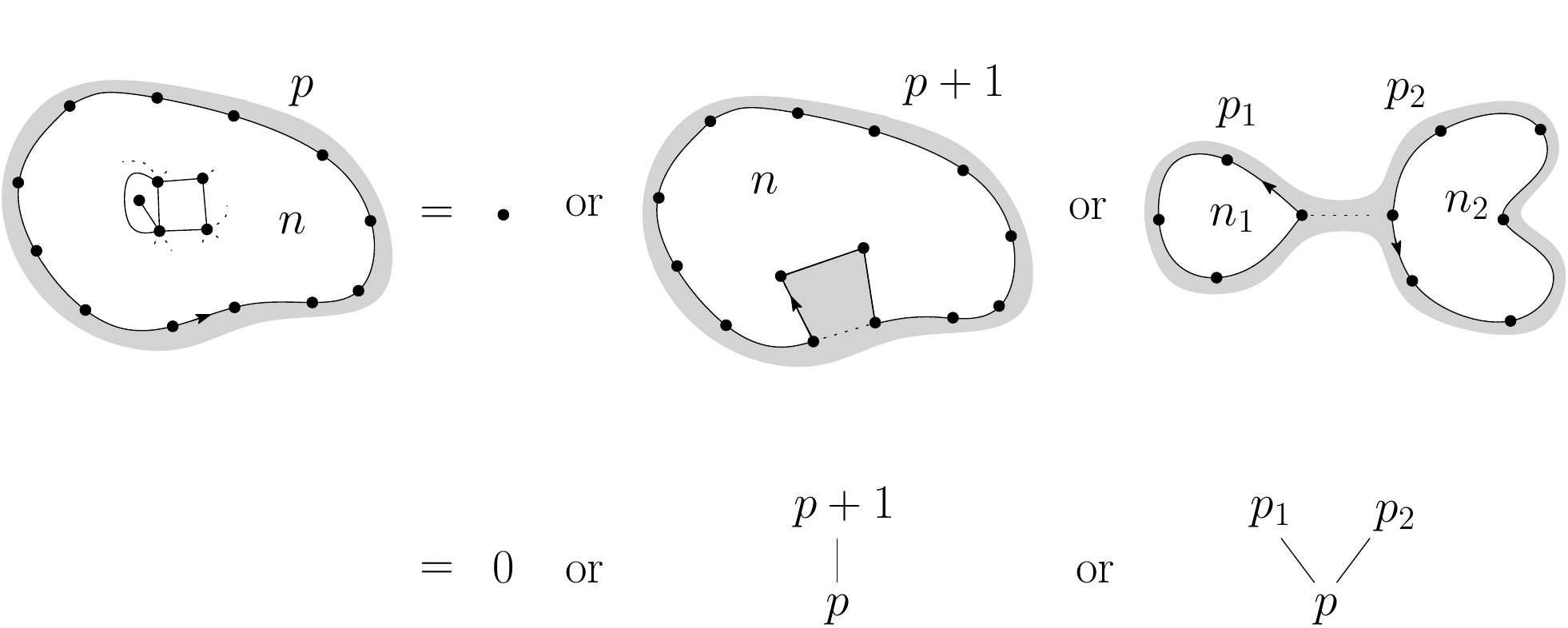}
 \caption{\label{fig:Tutte} Tutte's recursive decomposition for quadrangulations and bellow the local correspondence in the peeling tree. In the case when the removal of the root edge splits the quadrangulation (right), we put the component attached to the origin of the root edge on the left of the peeling tree.}
 \end{center}
 \end{figure}

\begin{figure}[!h]
 \begin{center}
 \includegraphics[width=14.5cm]{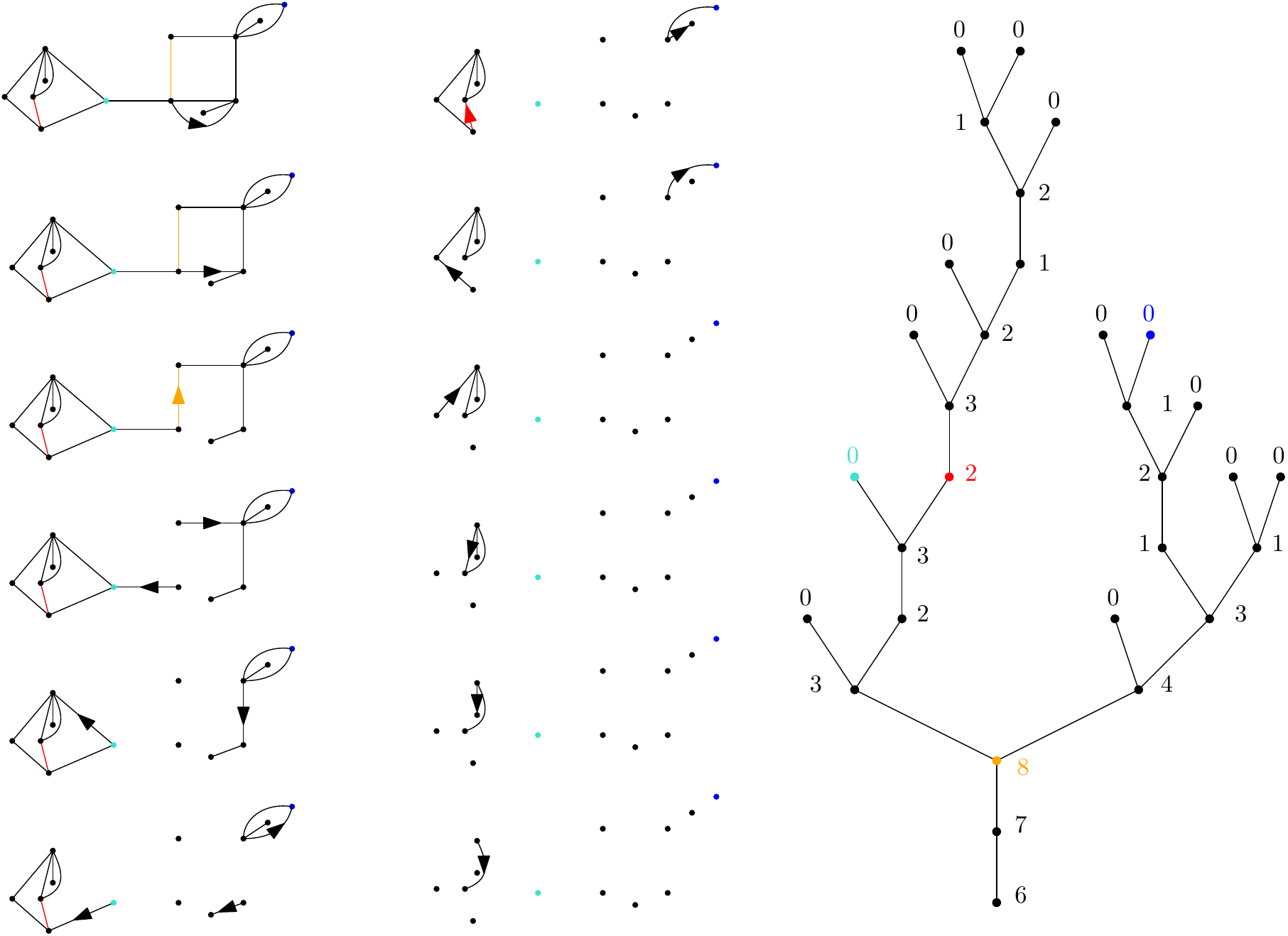}
 \caption{Step-by-step example of the peeling process on a quadrangulation with $n=12$ vertices and a boundary of degree $2 p$ with $p = 6$ and its corresponding peeling tree. In blue, two vertices of the quadrangulation and the two corresponding leaves with label $0$ in the peeling tree. In orange and rouge, two edges and their corresponding inner vertices in the peeling tree.\label{fig:peel-quad}}
 \end{center}
 \end{figure}

\paragraph{Enumeration of planar maps.}
In this paper, we will focus on \emph{planar maps}, which are finite connected graphs (possibly with loops and multiple edges) properly embedded in the two-dimensional sphere, seen up to continuous deformation. To avoid symmetries, all planar maps will be rooted at a distinguished oriented edge. The number of edges incident to a face i.e.\ that of its underlying polygon, is called the \emph{degree} (or sometimes the \emph{perimeter} or \emph{length}) of this face. 
We will particularly study the case of planar \emph{quadrangulations} and \emph{triangulations} which are maps where all faces are quadrangles (resp.\ triangles). Sometimes our maps will have a boundary i.e.\ a distinguished face which may have a different degree (but this degree has to be even in the case of quadrangulations), in which case the root edge lies on the boundary with the distinguished face lying to its right. This face will often be drawn and denoted as the external face.  \\
One possible way to enumerate quadrangulations or triangulations is to use Tutte's method based on a recursive decomposition which is obtained by removing the root edge. The point is that when removing an edge from a quadrangulation, it may not be a quadrangulation anymore but a quadrangulation with a boundary. 
More precisely, when erasing the root edge of a quadrangulation with a boundary with degree $2 p$ with $p \geq 1$, one of these two possible events occurs:
\begin{itemize}
\item either the quadrangulation stays connected, which means that one discovers a new face of the initial quadrangulation, and one re-roots the quadrangulation at the left-most edge of the new face.
\item or the deletion of the root edge disconnects the quadrangulation and one gets two quadrangulations with half-perimeter $p_1 \geq 0$ and $p_2 \geq 0$ such that $p_1+ p_2 = p-1$, that one can re-root  using the two endpoints of the removed edge, see Figure \ref{fig:Tutte}.
\end{itemize}
Note that we considered here the half-perimeter since the quadrangulations are bipartite and we only encounter even boundaries during the exploration.
In the case $p = 0$, the map is just a vertex map i.e.\ the map composed of one single vertex and no edge. One can iterate Tutte's decomposition in each quadrangulation with a boundary until one obtains a collection of vertex maps.
This edge-by-edge peeling exploration of a quadrangulation can be encoded in a so-called \emph{peeling tree}\footnote{In random map theory, we often consider different ways to re-root the components, or \emph{peeling algorithms}, which may yield to different peeling trees, see \cite{CurStFlour}. For the sake of simplicity, we shall stick here to the rules presented above and do not use any other peeling algorithm.} by recursively labeling the vertices with the half-perimeters of the boudaries, see Figure \ref{fig:Tutte}. Conversely, given a labeled plane tree whose labeling follows the appropriate rules described above, one can recover the initial quadrangulation. Note that the vertices of the peeling tree with label $0$ are the leaves and correspond to the vertices in the quadrangulation, whereas the other vertices correspond to the edges of the map, see Figure \ref{fig:peel-quad}.

Tutte then obtained recursive equations for the number of quadrangulations using the boundary length as a catalytic variable which can be summarize into the following equation on $ \mathbf{Q}$ the bivariate generating function of quadrangulations with a boundary where the variable $x$ counts the number of vertices and $y$ counts the half-perimeter of the boundary
\begin{equation}\label{eq:tutte-quad} \mathbf{Q} =x+ y \mathbf{Q}^2 +  \frac{1}{y} \left( \mathbf{Q} -x- y \mathfrak{Q}\right),
\end{equation}
where $\mathfrak{Q} = [y^1] \mathbf{Q}$ is the (univariate) generating function of quadrangulations with a boundary of length~$2$. On the right, the term $x$ stands for the vertex max, and the term $y \mathbf{Q}^2$ encompasses the case when the removal of the root edge splits the quadrangulation. The remaining term corresponds to the case when the quadrangulation stays connected in which case the new quadrangulation has at least a boundary of length $4$. This equation characterizes the power series $\mathbf{Q}$ and has been solved explicitly using the so-called ``quadratic method" and then its generalization the kernel method introduced by Bousquet-Mélou and Jehanne~\cite{BMJ06}. 
 
Since then, other methods have been developed to enumerate maps: via matrix integrals \cite{BIPZ78,tHo74}, bijections with other labels trees ``\`a la Schaeffer" \cite{CV81,schaeffer1997bijective} or correspondence with the KP hierarchy \cite{GJ08}. Our work concentrates on peeling trees but uses another method to enumerate them which is based on Lackner--Panholzer last car decomposition of parking trees \cite{LaP16}. This link between parking models and maps was already suggested by Panholzer in \cite[Remark 2]{panholzer2020combinatorial} and by the author and Curien in \cite[Section 8]{contat2021parking}. Panholzer found remarkable explicit enumeration formulas of parked trees for a large class of combinatorial models, one of which is linked to the enumeration of non-decomposable maps. The non-decomposable maps also have links with description trees \cite{CoriSchaefferDescription} of Cori and Schaeffer whose construction shares similarities with that of parking trees which we explain now. The link between bipartite planar maps and ``degree trees" pointed out by Fang in \cite{fang2021bijective} also supports this strong link between parking trees and map models.

\paragraph{Parking on trees.} Our decomposition uses an idea introduced by Lackner and Panholzer in \cite{LaP16} in the context of parking trees. Let us recall this model for the readers' convenience although we shall not use it in this paper. Let $ \mathbf{t}$ be a finite (plane) rooted tree which is our parking lot. Each vertex represents a parking spot which can accommodate at most one car. We let then cars arrive on the vertices of $ \mathbf{t}$. Each car tries to park on its arrival node but if the spot is already occupied, it drives towards the root and parks as soon as possible. An important property of this model is its Abelian property: the final configuration, the flux of cars which go through a given edge and the outgoing flux of cars do not depend upon the order chosen to park the cars. In particular, on can recover the initial configuration of cars from the final configuration of parked cars and flux on the edges.   This model undergoes a phase transition and was first studied on a directed line \cite{konheim1966occupancy} and raised recently a growing interest, especially on random tree models with an increasing level of generality \cite{LaP16,GP19,chen2021parking,CH19,contat2020sharpness}.

\begin{figure}[!h]
 \begin{center}
 \includegraphics[width=14cm]{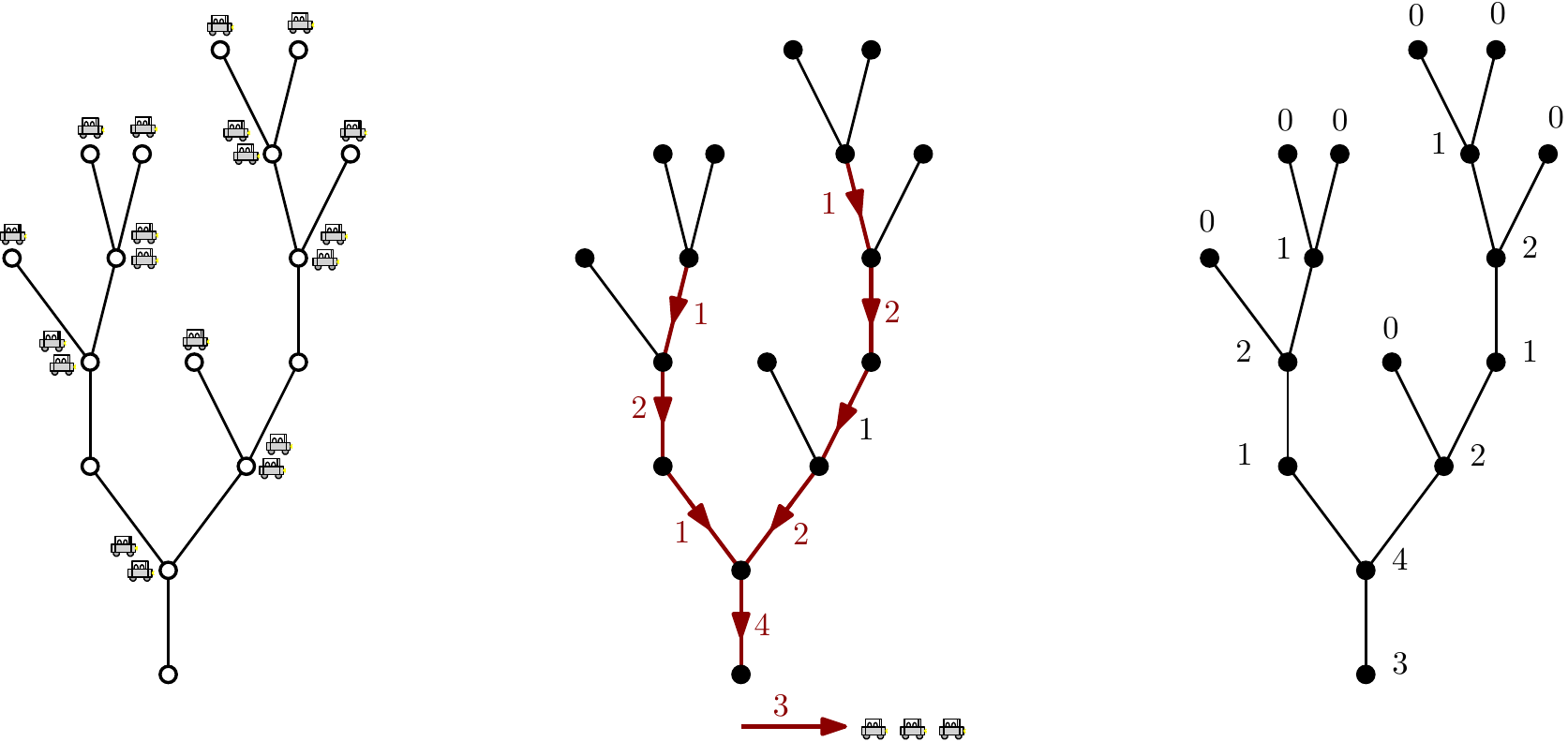}
 \caption{\label{fig:parking} Left, an example of a fully parked tree with 16 vertices and 19 cars, three of which can not park. We put exactly one car arriving on every leaf, two cars on every vertex with two children, and no car on vertices with one child. In the middle, the corresponding final configuration where the edges are decorated with the flux of cars (when non-zero). When pushing up the labels to the above vertices, one obtains a labeled tree which follows the same local rules as the peeling tree of quadrangulations (right).}
 \end{center}
 \end{figure}

We will focus here on the connected components of parked vertices in the final configuration in this model, which are called \emph{fully parked trees}. They consist of plane trees with a decoration of cars so that all vertices accommodate a car and some cars possibly contribute to the outgoing flux. We can decorate the edges of those trees with the flux of cars and in a specific case of trees and car arrivals, the labeled trees obtained that way (pushing up the labeling from the edges to the vertices above) are very similar to our peeling trees since the rules for the labels are identical, see Figure \ref{fig:parking}. Specifically consider a fully parked tree where the vertices have $0$, $1$ or $2$ children and such that exactly one car arrives on each leaf (vertex with $0$ children), no car arrives on vertices with one child and two cars arrive on vertices with two children. Then the leaves will all get label $0$ since one car arrives and parks on each leaf and no car can come from above. If a vertex has one child with label $ \ell \geq 1$ (and no car arriving on it), then one of the $ \ell $ cars arriving from above parks and the vertex will get label $ \ell - 1$. Lastly, a vertex with two children with label $ \ell_1$ and $ \ell_2 \geq 0$ has two cars arriving on it and $\ell_1 + \ell_2$ coming from above. One of them parks and it remains $ \ell_1 + \ell_2+1 $ cars contributing to the flux on the edge below. Those local rules are exactly the same as that of peeling trees of quadrangulations and we can thus match bijectively these fully parked trees (with prescribed car arrivals) with quadrangulations with a boundary.
 
To enumerate such plane general trees, Chen \cite{chen2021enumeration} uses a method which is similar to Tutte's since he decomposes the fully parked trees at their root using the outgoing flux of cars as a catalytic variable. 

The technique introduced by Lackner \& Panholzer \cite{LaP16} and deepened by Panholzer in \cite{panholzer2020combinatorial} is different: it consists in a decomposition of the initial fully parked tree according to the parking spot of  a distinguished car seen as the ``last" car. Thanks to the Abelian property, we can imagine that we first park all the cars but this distinguished car and that it arrives at last and parks on a vertex, which was empty before the last car arrived. And when removing this last car, the trees which are attached to this free spot are also fully parked trees. See also \cite[Section 8.2]{contat2021parking} for the decomposition of fully parked trees with a different notion of components and with possibly a flux of outgoing cars.

The heart of this work is to adapt the ``last car decomposition" to the enumeration of planar maps, more precisely of the peeling trees of planar maps. This new decomposition is explained in Section \ref{sec:decomp} in the case of quadrangulations and in Section \ref{sec:triang} for triangulations with zero, one or two boundaries. In both cases, this decomposition enables us to establish new equations on the corresponding generating functions which we state now.

\paragraph{Quadrangulations.}
We introduce $Q_{n}$ be the number of rooted quadrangulations (without boundary) with $n$ vertices, and we denote by $ \mathfrak{Q}$ its corresponding generating function $$ \mathfrak{Q} (x) := \sum_{n \geq 2}Q_{n} x^n = x^2 +2 x^3 + 9x^4 + 54x^5+378x^6 + \cdots.$$ 
We already encountered $ \mathfrak{Q}$ as the generating function of quadrangulations with a boundary of perimeter $2$. Indeed we will see later  that we can transform the root edge of a quadrangulation without boundary to get (bijectively) a quadrangulation with boundary of degree $2$.
By convention, the map with two vertices linked by an edge is considered as a quadrangulation (without boundary) with $0$ face, which explains the term $x^2$. 
\begin{theorem}\label{thm:quad0}
Writing $\mathfrak{Q}^{\bullet} = x \mathfrak{Q}'(x)$, the ``last car decomposition" translates into the following equation
\begin{equation}\label{eq:quad} \mathfrak{Q}^{ \bullet} = 2x^2 + 6x \left(\frac{ \mathfrak{Q}^{\bullet}-  \mathfrak{Q}}{1 - \mathfrak{Q}^{\bullet}/x}\right),  \end{equation}
which characterizes $  \mathfrak{Q}$ and which is equivalent to the following recursive equation:  $Q_{2} = 1$ and for $n \geq 3$,
\begin{equation} \label{eq:rec-quad}n Q_{n}=  \sum_{k = 2}^{n-1} k (n+1-k) Q_{k}  Q_{n+1-k}+ (4n-10) Q_{n-1}.\end{equation}
\end{theorem}
Notice that Equation \eqref{eq:quad} or equivalently \eqref{eq:rec-quad} does not make use of quadrangulations with a boundary. They are fundamentally different from those of Tutte \eqref{eq:tutte-quad} which necessitated the introduction of a \emph{catalytic variable} $y$ to write the equation on $\mathbf{Q}$ and eventually characterize $ \mathfrak{Q}$ alone. 
The form of the above equations may remind the knowledgeable reader the decomposition obtained via the KP hierarchy, see for example \cite[Corollary 1]{louf2019new}. But it seems that those equations can not be deduced one from another. Indeed with our notation, the equation coming for the KP hierarchy is 
\begin{equation} \mathfrak{Q}^{\bullet} - \mathfrak{Q} = 4 x ( 2  \mathfrak{Q}^{ \bullet} -3 \mathfrak{Q}) + 3 (2  \mathfrak{Q}^{\bullet} - 3 \mathfrak{Q})^2 +x^2.
\end{equation}
Our decomposition also allows us to deduce a recursive decomposition of quadrangulations with a boundary. Let $Q_{n}^{(p)}$ be the number of rooted quadrangulations with a boundary of length $2 p$ and $n$ vertices, and we denote by $ \mathbf{Q}$ its corresponding bivariate generating function $$ \mathbf{Q} (x,y) := \sum_{n \geq 1}\sum_{ p = 0}^{n} Q_{n}^{(p)} x^n y^{p} = x + y \left( x^2 +2 x^3 + 9x^4 +\cdots \right) + y^2 \left( 2 x^3 + 9x^4 + 54x^5 +\cdots \right)  + \cdots.$$ 
By convention, we consider the single isolated vertex as a planar quadrangulation with a boundary of perimeter $0$ (and $0$ face), which explains the term $x$. 
The last car decomposition gives the following differential equation on $ \mathbf{Q}$
\begin{equation}\label{eq:edquad} \mathbf{Q}^{\bullet} = x + 6y \mathbf{Q}^{\bullet} \left(\frac{ \mathfrak{Q}^{\bullet}- \mathfrak{Q}}{1 - \mathfrak{Q}^{\bullet}/x}\right) + 2xy\left(3\mathbf{Q}^{\bullet}- 2\mathbf{Q}- y\partial_{y}\mathbf{Q}\right), \end{equation}	
where $ \mathbf{Q}^{\bullet} = x \partial_{x}\mathbf{Q}$. This equation equation is not as ``simple" as Tutte's Equation \eqref{eq:tutte-quad} since it also involves the partial derivatives of the generating function $ \mathbf{Q}$.

Lastly we can also consider quadrangulations with multiple boundaries that are enumerated by Tutte's slicing formula, see \cite[Theorem 3.4]{CurStFlour}. We will only deal with the case of two boundaries at the end of Section \ref{sec:decomp} but the general case can be obtained by stacking  and gluing the appropriate number of peeling trees with the same procedure. 

\paragraph{Triangulations.}
An adaptation of the ``last car" technique gives similar results in the case of triangulations. As above we start by the case without boundary, and the map with two vertices linked by an edge is considered as a triangulation with $0$ face. We denote by $T_n$ the number of rooted triangulations with $n$ vertices without boundary and we denote its generating series~by $$\mathfrak{T} (x) = \sum_{n \geq 2} T_n x^n = x^2 + 4 x^3+ 32x^4+ 336x^5 + \dots$$
\begin{theorem}The last car decomposition yields the following equation
\begin{equation} \label{eq:tri}
3\mathfrak{T}^{\bullet} - 4 \mathfrak{T} = 2\left(\frac{ \mathfrak{T}^{\bullet}- \mathfrak{T}}{1 - \mathfrak{T}^{ \bullet}/x}\right),
\end{equation}
where $ \mathfrak{T}^{\bullet} = x \mathfrak{T}'(x)$. Assuming $T_2 = 1$, this equation characterizes~$\mathfrak{T}$ and is equivalent to the recursive equation which is similar to \eqref{eq:rec-quad}
\begin{equation}\label{eq:tri-rec} T_{n}:= \frac{1}{n-2} \sum_{k=2}^{n-1} (3k-4)(n+1-k)T_kT_{n+1-k}.
\end{equation}
\end{theorem}
As above, this equation may remind the reader the following equation coming from the KP hierarchy:
\begin{equation} \mathfrak{T}^{\bullet} - \mathfrak{T} = (6 \mathfrak{T}^{\bullet} - 8 \mathfrak{T} +x)^2,
\end{equation}
 see \cite[Theorem 5.4 and Equation 45]{GJ08}. However we were not able to deduce one from another by simple computations. 
We now let $T_{n}^{(p)}$ be the number of rooted triangulations with a boundary of length $ p$ and $n$ vertices in total, and we denote by $ \mathbf{T}$ its corresponding generating function $$ \mathbf{T} (x,y) = \sum_{n \geq 1, p \geq 0} T_{n}^{(p)} x^n y^{p} = x + y \mathfrak{T}(x) +  \dots.$$
Then our last car decomposition shows that $ \mathbf{T}$ satisfies the following differential equation
\begin{equation} \label{eq:tri1bord}
6 \mathbf{T}^{\bullet} - 2 y \partial_{y}\mathbf{T} - 6\mathbf{T} + y \partial_{y} (y \mathbf{T}) =  \frac{4y}{x} \mathbf{T}^{\bullet} \left(\frac{ \mathfrak{T}^{\bullet}- \mathfrak{T}}{1 - \mathfrak{T}^{\bullet}/x}\right) + y \left(4\mathbf{T}^{\bullet}- 3\mathbf{T}- y\partial_{y}\mathbf{T}\right). \end{equation}
where $ \mathbf{T}^{\bullet} = x \partial_{x}\mathbf{T}$, which is the analogue of Equation \eqref{eq:edquad}. We will see later that a small local transformation on the triangulations shows that $ \mathfrak{T} = [y^1]\mathbf{T}=  [y^2]\mathbf{T}= [y^3]\mathbf{T}$. 
We can also adapt our decomposition in the case of triangulations with two boundaries, see the end of Section \ref{sec:triang}.

\paragraph{Acknowledgements.} We acknowledge support from ERC 740943 \textit{GeoBrown}. We thank Thomas Budzinski and Baptiste Louf for stimulating discussions. I am also grateful to Nicolas Curien for his precious suggestions.

\section{``Last Car" decomposition of quadrangulations}\label{sec:decomp}

We have given above an insight into the peeling exploration technique on quadrangulation with a boundary.  Let us precisely explain the correspondence between quadrangulations with a boundary and their peeling tree. Recall that the removal of the root edge in a quadrangulation can produce two possible events: either it stays connected and the half-perimeter rises by $1$, or it splits in two parts and the sum of the two half-perimeter is prescribed, see Figure \ref{fig:Tutte}.  \\

\noindent\fbox{\parbox{14cm}{\textbf{Peeling tree of quadrangulations with one boundary.}
 To each rooted quadrangulation with $n$ vertices and with a boundary of perimeter $2 p$, we can match bijectively its peeling tree, which is a labeled plane tree whose vertices have between zero and two children including $n$ leaves. The labeling satisfies the following rules:
\begin{itemize}
\item The root vertex has label $p$.
\item The leaves have label $0$.
\end{itemize} 
\vspace{-0.28cm}
\begin{itemize}[label=$\star$]
\item  \hypertarget{itm:loc-quad}{\emph{(local rules)}} All inner vertices have a label $\ell \geq 1$ and either a child with label $\ell+1$ or two children with labels $ \ell_1, \ell_2 \geq 0$ such that $\ell_1+\ell_2 -1= \ell$. See Figure \ref{fig:peel-quad}.
\end{itemize} 
Note that the labels of the vertices can be seen as the half-perimeter of the successive boundaries during the exploration, and that the local rules imply that it has $2n-p-2$ inner vertices, which is also the number of edges of the quadrangulation. 
}} \\

\begin{figure}[!h]
 \begin{center}
 \includegraphics[width=12.5cm]{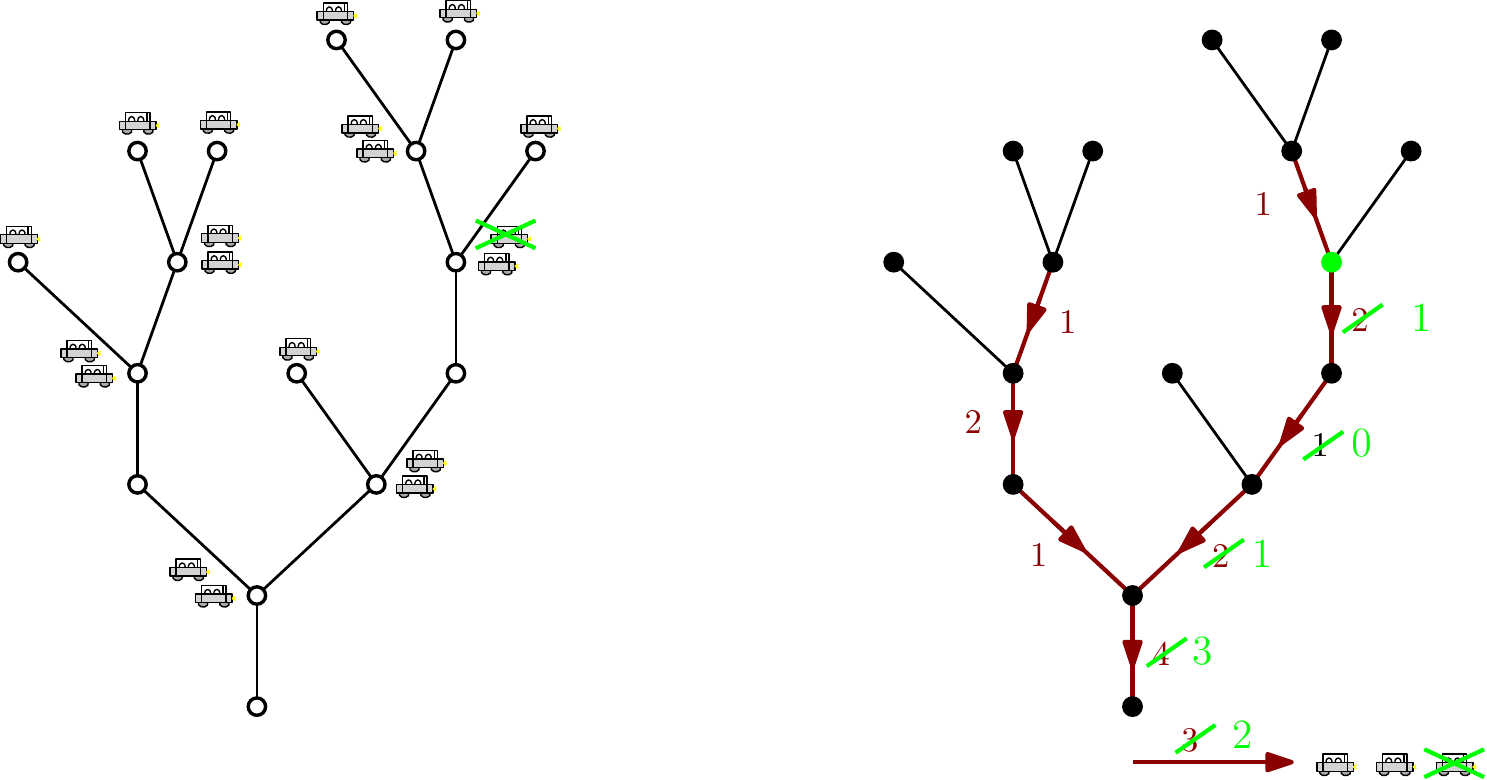}
 \caption{\label{fig:parklast} The effect of removing a car (seen as the last car) in a fully parked tree. We remove here the car in green on the left and imagine it parked last. On the right, it decreases by $1$ the flux of the edges between the green point (location of the car arrival) and the root edge.}
 \end{center}
 \end{figure}

We will now concentrate on such peeling trees and (almost) forget about their interpretation in terms of maps. We want to apply a ``last car" decomposition to such a tree. Let us first explain the influence of removing a car from a parking tree or fully parked tree.
Recall that a fully parked tree is a rooted tree together with a car configuration where all parking spots are occupied (eventually with an outgoing flux) and that we can label each vertex by the flux of cars that go through the edge just below (or the outgoing flux for the root vertex). Imagine now that one removes a distinguished car from this tree, or more precisely that we first park all the cars but this one and then try to park this distinguished car so that we can easily remove it. Take for example a car which contributed to the outgoing flux. Then the effect of this removal is that the flux in the edges of the branch between the arriving spot of this car and the root vertex has decreased by $1$, see Figure \ref{fig:parklast}.

In the next sections, we will try to subtract $1$ to the labels in a branch of a peeling tree. We will see that in the inner vertices of the branch, this transformation preserves the local rules. 

\subsection{Quadrangulations without boundary or with a boundary of length $2$.}
\begin{figure}[!h]
 \begin{center}
 \includegraphics[width=10cm]{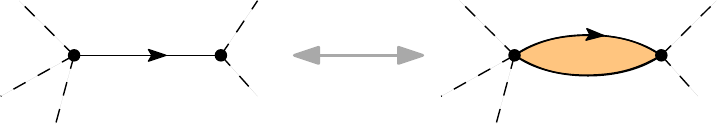}
 \caption{\label{fig:root-transf-quad} Root transform for quadrangulations (or more generally for bipartite maps). On the left the initial root edge and on the right the distinguished face of degree $2$ together with the new root edge.}
 \end{center}
 \end{figure}

To enumerate quadrangulations without boundary we actually first transform them into quadrangulations with a boundary of degree $2$. To this end, the standard trick is to cut along the root edge and ``open" it, see Figure \ref{fig:root-transf-quad}. This transformation does not affect the number of vertices of the initial quadrangulation and is a bijection between quadrangulations without boundary and with a boundary of perimeter $2$ with the same number of vertices. Since we only consider planar maps, Euler's formula implies that a quadrangulation without boundary has $2n-4$ edges (and $n-2$ faces), and applying the root transform only increases the number of edges and faces by $1$.  

\begin{figure}[!h]
 \begin{center}
 \includegraphics[width=14cm]{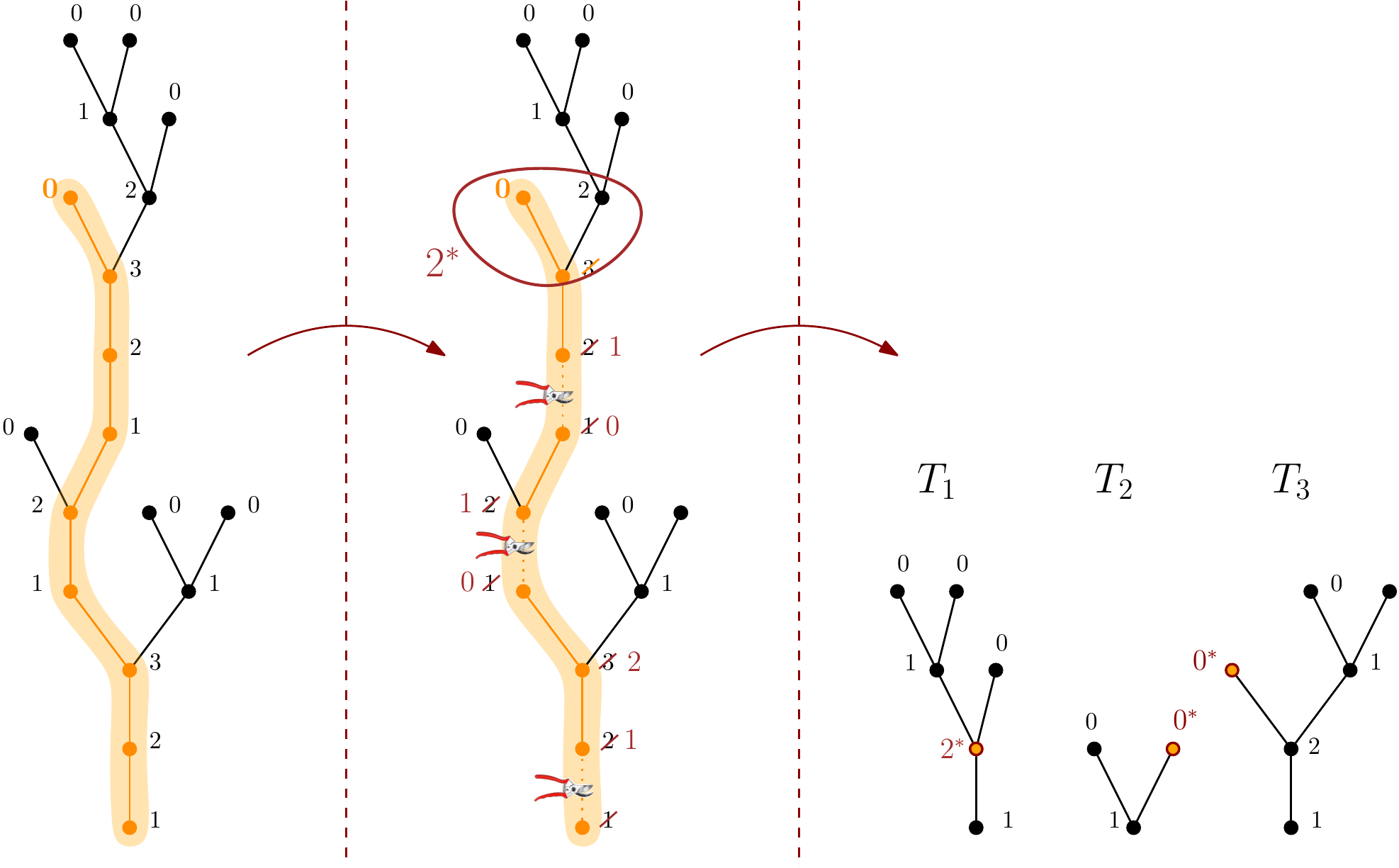}
 \caption{\label{fig:transfoQuad2}Example of the last car decomposition on a peeling tree of quadrangulation with $n=7$ leaves into $k=3$ trees with $n+k-2 = 8$ leaves in total. On the left, the initial peeling tree where the distinguished leaf is displayed in orange as well as the branch between this leaf and the root. In the middle, we point out the needed transformations when removing one in this branch to preserve the local rules: the father and sibling of the distinguished leaf are contracted into a (marked) vertex labeled $2$, we remove the initial root (or cut just above), and we cut above the two vertices which get label $0$. On the right, the resulting sequence of three trees.}
 \end{center}
 \end{figure}

Thanks to this trick, we now only need to concentrate on peeling trees of quadrangulations with a boundary of length $2$ i.e. those whose root has label $1$. \\
Inspired by removing the last car in fully parked trees \cite{LaP16}, we want to remove $1$ to all labels in a ``branch" of the peeling tree. To do this, we need to distinguish a leaf of the tree, which corresponds to a vertex in the initial quadrangulation and we consider the branch between this leaf and the root vertex. The key observation is that subtracting $1$ to all labels of this branch preserves the local rules of the tree which we described in (\hyperlink{itm:loc-quad}{$\star$}), but two issues may appear. 
\begin{itemize}
\item First, the initial distinguished leaf cannot get a label $-1$, but it has a parent with some label $k \geq 1$ and a sibling with label $k-1$, so that we can just contract  these three vertices into a vertex with label $k-1$. We distinguish it to remember where we removed the distinguished leaf.
\item The other possible issue is that the vertices with label $1$ will get label $0$ and therefore have to be leaves. This is also easily solved by cutting the edges just above the new $0$'s and distinguishing them to remind the location of the cuts, see Figure \ref{fig:transfoQuad2}. 
\end{itemize}

Lastly, when $n \geq 3$, the root vertex of the tree has always a single child with label $2$ so that we just remove the initial root and root the new tree at this child which gets label $1$. \\
To summarize, our transformation converts a peeling tree starting from $1$ with $n \geq 3$ leaves, marked at one leaf, into a sequence of $k$ trees  for some $k \geq 1$, with root label $1$, with $n+k-2$ leaves in total and such that the first tree has a distinguished vertex (leaf or inner vertex) and the eventual following trees have a marked leaf. 
Conversely, given such a sequence, we can recover the initial  peeling tree by gluing successively the root of a tree with the distinguished leaf of the next tree of the sequence and add $1$ to the labels in the appropriate branch.
Noting that a quadrangulation with a boundary of perimeter $2$ and $n$ vertices has $2n-1 - 2$ edges hence its peeling tree has $3n-3$ vertices in total, we obtain the following differential equation on $ \mathfrak{Q}$:
\begin{equation*} \mathfrak{Q}^{ \bullet} = 2x^2 + 2x \left(\frac{ 3\mathfrak{Q}^{\bullet}- 3 \mathfrak{Q}}{1 - \mathfrak{Q}^{\bullet}/x}\right),  \end{equation*}
which is Equation \eqref{eq:quad}. Let us explain this in more details:
\begin{itemize}
\item The term $2x^2$ stands for the map with two vertices connected by the root edge, which is a quadrangulation by convention.
\item The factor $2x$ comes from the fact that there are two possible ways to ``reglue" the initial marked leaf with label $0$ on the first tree of the transformed sequence (left or right). 
\item The factor $3\mathfrak{Q}^{\bullet}- 3 \mathfrak{Q}$ comes from the first tree in the sequence, which has a distinguished vertex, and is divided by $1 - \mathfrak{Q}^{\bullet}/x$ for all possible cuts. 
\end{itemize}

By identifying the coefficient in the above equation, we obtain Equation \eqref{eq:rec-quad} which concludes the proof of Theorem \ref{thm:quad0}. This is a new equation which characterizes A000168 in Sloane online encyclopedia for integer sequences.

\subsection{Quadrangulations with a boundary of degree $ 2 p \geq 4$ }

The transformation which we described above also works for quadrangulations with a boundary of perimeter $2 p \geq 4$. The only difference is that in that case, we do not remove the initial root of the tree.  Starting from a peeling tree with $n \geq 3$ leaves, a root labeled $p\geq 2$ and a distinguished leaf, we apply our last car decomposition in the branch between the distinguished leaf and the root of the tree. We obtain a sequence of $k \geq 1$ peeling trees for some $k \geq 1$ such that: 
\begin{itemize}
\item the first tree is marked at a vertex (leaf or not),
\item all eventual following trees have a marked leaf,
\item all but the last tree have a root label $1$ and the last tree has a root labeled $p - 1 \geq 1$,
\item the trees have $n+k -2$ leaves with label $0$ in total.
\end{itemize}

\begin{figure}[!h]
 \begin{center}
 \includegraphics[width=14cm]{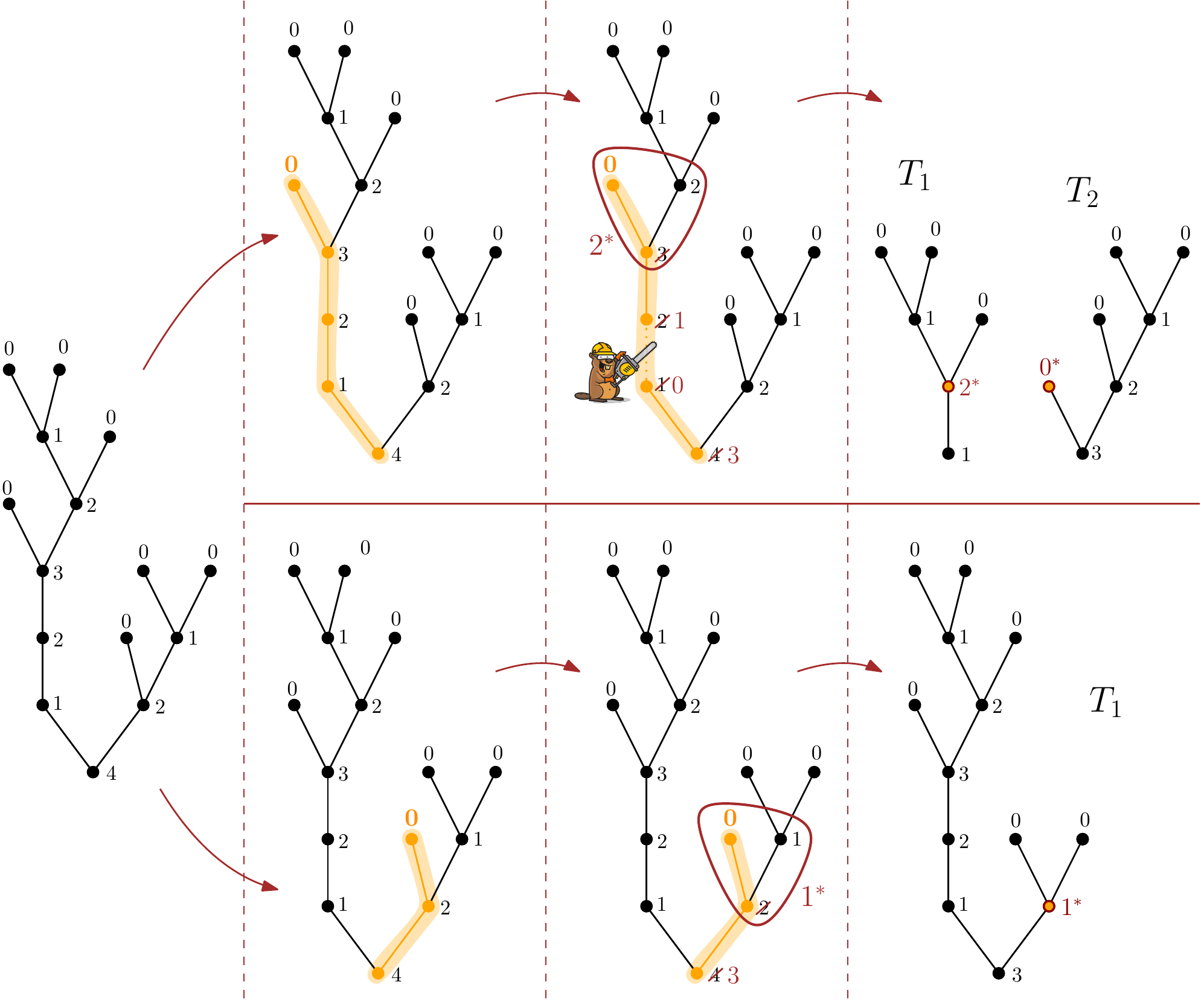}
 \caption{\label{fig:quad-bord} Two examples of our transformation on the same tree with root labeled $p=4$ and $n=7$ leaves with two different distinguished leaves. On the top, the transformation produces $k=2$ trees and have $n-k-2=7$ leaves in total. The first one has a distinguished vertex and root labeled $1$, and the second one has a distinguished leaf and root labeled $p-1=3$. At the bottom, we get $k=1$ tree with $n-k-1=6$ leaves, a distinguished vertex and a root labeled $p-1=3$. In particular, the first and last tree are the same. }
 \end{center}
 \end{figure}

See Figure \ref{fig:quad-bord}. Note that the first and the last tree can be confounded.
We then get the differential Equation \eqref{eq:edquad} on $ \mathbf{Q}$ that we recall here:
\begin{equation*}\mathbf{Q}^{\bullet} = x + 6y \mathbf{Q}^{\bullet} \left(\frac{ \mathfrak{Q}^{\bullet}- \mathfrak{Q}}{1 - \mathfrak{Q}^{\bullet}/x}\right) + 2xy\left(3\mathbf{Q}^{\bullet}- 2\mathbf{Q}- y\partial_{y}\mathbf{Q}\right), \end{equation*}	
where $ \mathbf{Q}^{\bullet} = x \partial_{x} \mathbf{Q}$ and $ \mathfrak{Q} = [y^1] \mathbf{Q}$, (resp.\ $ \mathfrak{Q}^{\bullet} = [y^1] \mathbf{Q}^{\bullet}$). Indeed, 
\begin{itemize}
\item The term $x$ corresponds to the map composed of a single vertex and no edge whose peeling tree is simply the tree with one vertex labeled $0$.
 \item The second term corresponds to the case when our transformation gives more than one tree (and in particular, the first and the last trees are not the same), including the case of quadrangulations with a boundary of perimeter $2$ where the last tree after transformation is the tree with one vertex labeled $0$. In that case, the last tree is a tree with root label $p-1$ and a marked leaf, which corresponds to the factor $y \mathbf{Q}^{ \bullet}$. The first tree has a distinguished vertex and its root has label $1$, hence the factor $ 3  \mathfrak{Q}^{ \bullet}- 3 \mathfrak{Q}$ as in the case of a boundary of perimeter~$2$. The factor $1/(1 - \mathfrak{Q}^{\bullet}/x)$ stands for the other possible cuts. There is also a factor $2$ which comes from the fact that there are two possible ways to reglue the initial distinguished leaf labeled $0$ on the first tree.
\item The last term corresponds to the quadrangulations with $n \geq 2$ vertices where we get only one tree when applying our transformation on its peeling tree. In that case the resulting tree has $n-1$ leaves, a root labeled $p-1$ hence $2(n-1)-( p - 1)-2 = $ inner vertices (edges in the quadrangulation) and $ 3(n-1) - 2 - ( p -1)$ vertices in total, one of them is distinguished explaining the factor $3\mathbf{Q}^{\bullet}- 2\mathbf{Q}- y\partial_{y}\mathbf{Q}$.
\end{itemize}

\subsection{Quadrangulations with two boundaries}
\begin{figure}[!h]
 \begin{center}
 \includegraphics[width=15cm]{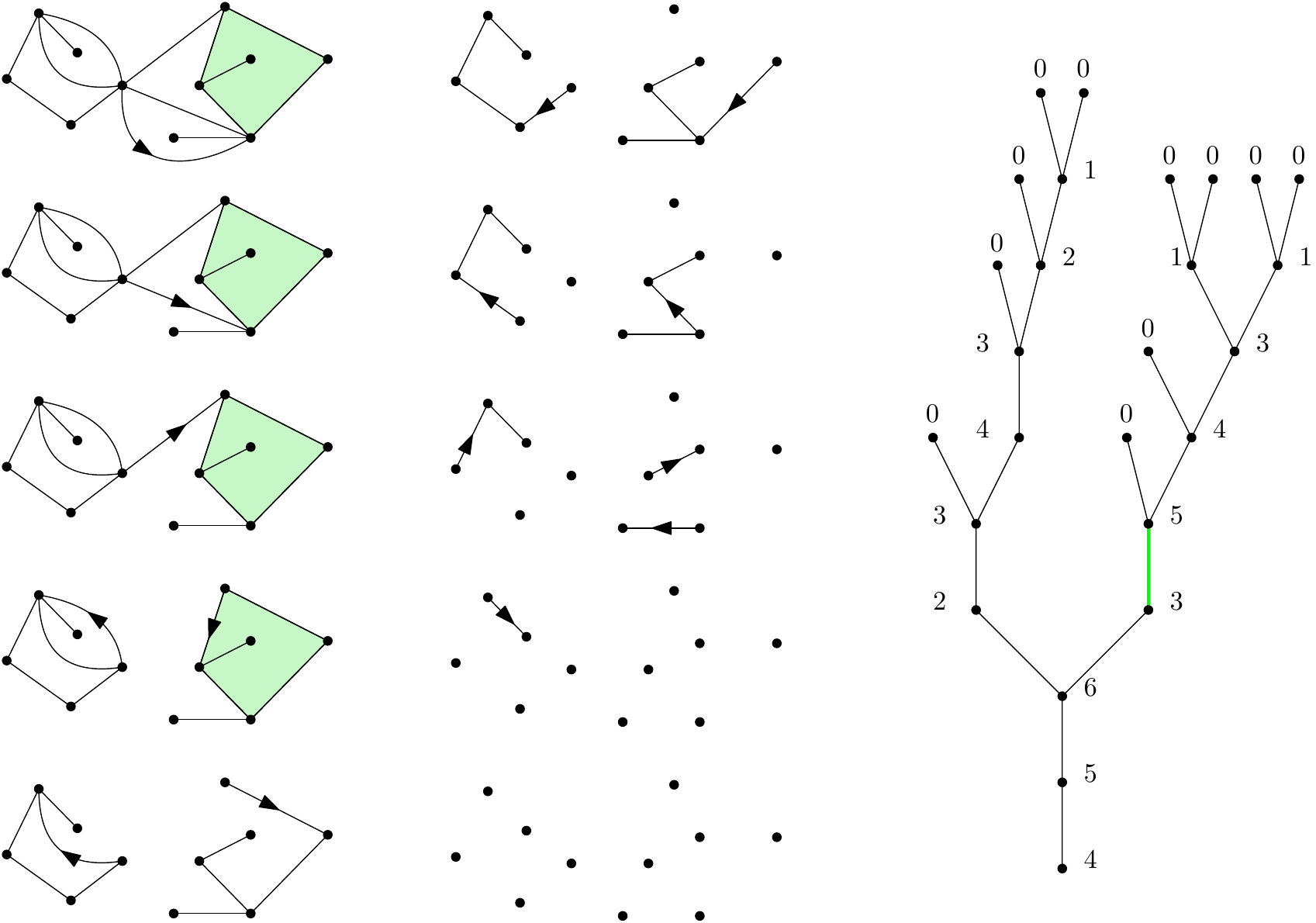}
 \caption{\label{fig:quad2bord}Step-by-step peeling exploration of a quadrangulation with a boundary of perimeter $8= 2\cdot 4$ and a distinguished hexagon in green. The discovery of this face leads to a transition from a boundary of half-perimeter $3$ to half-perimeter $3+3-1= 5$ (green edge in the peeling tree on the right). Removing this green edge in the peeling tree on the right, the above part is a usual peeling tree of a quadrangulation with a boundary of half-perimeter $5$ and $6$ leaves. The bottom part starts from a label $4$, has $5$ leaves with label $0$ and a distinguished leaf with label $3$.}
 \end{center}
 \end{figure}

We can also give a recursive equation for the number of quadrangulations with two boundaries. 
We interpret the first boundary as above but the second boundary as a distinguished face of perimeter $2q$, and denote by ${Q}_{n}^{(p,q)}$ the number of rooted quadrangulations with a boundary of length $2 p$ and a distinguished face of perimeter $2q$ and $n$ vertices in total. 
In the peeling exploration, the discovery of the distinguished face corresponds to the event of seeing a boundary of half-perimeter $r \geq 1$ which becomes a boundary of half-perimeter $r+q-1$ when removing an edge.
This matches in the peeling tree to a vertex with label $ r \geq 1$ which gives birth to a vertex labeled $r+q-1$ when discovering this distinguished face. The rest of the local transitions are exactly the same as before, see Figure \ref{fig:quad2bord}. Thus, we introduce ${A}_n^{(p,r)}$ the number of peeling trees of ``quadrangulations" with $n$ leaves labeled $0$, a distinguished leaf with label $r$ and a root labeled $p$. 

Indeed, if we decompose the peeling tree of a quadrangulation with two boundaries of length $2p$ and $2q$ according to this different transition, then for some $r \geq 1$, the above part is the usual peeling tree of a quadrangulation with one boundary of length $r+q-1$ and the other one has root label $p$ and a distinguished leaf with label $r$. And the total number of leaves with label $0$ should be $n$ if the quadrangulation had $n$ vertices. 
With this decomposition, we obtain 
\begin{equation} \label{eq:rec2bords}{Q}_{n}^{(p,q)} =  \sum_{r=1}^{n}\sum_{k = 0}^{n-r-q} {A}_{k}^{(p,r)} Q_{n-k}^{(r+q-1)}.
\end{equation}
Knowing the ${A}_n^{(p,q)}$'s, this equation allows us to deduce the $ {Q}_{n}^{(p,q)} $. Note that the sum on $r$ when $k=0$  encompasses the case where the bottom part has only one leaf with label $r$  and is just a straight vertical line (thus $ {A}_{0}^{(p,r)}$ equals $1$ if $r \geq p$ and $0$ otherwise).  

It only remains to enumerate peeling trees with a distinguished leaf with label $r \geq1$. For this, we can use our last car decomposition and subtract $1$ along the branch to the distinguished leaf. In that case, we do not need to transform the father and the sibling of the distinguished leaf. When $r=1$, then removing $1$ in the branch between this leaf and the root of the tree creates a leaf with label $0$ and thus a usual peeling tree of a quadrangulation with one boundary. Hence we set ${A}_{n}^{(p,0)} = (n+1) Q_{n+1}^{(p)}$  when $n \geq\max(1,p-1)$. \\
When $r \geq 2$, the distinguished leaf just get label $r-1 \geq 1$. 
We then obtain the following recursive equation:
\begin{eqnarray}
\forall r \geq 1, \forall p \geq 1, \forall n \geq 1, &{A}_{n}^{(p,r)} &= {A}_{n}^{( \max(p-1,1) , r-1)} + \sum_{k=0}^{n-1} {A}_{k}^{(1,r-1)} {A}_{n-k}^{(p,1)}.\label{eq:recTreeQuadR}
\end{eqnarray}
The first term on the right corresponds to the case when there is no $1$ in the branch so that there is no cut in the decomposition. When $p=1$, we remove the initial root whereas we do not remove it when $p \geq 2$ so that the root of the new tree has label $p-1$, which explained the index $ \max(p-1,1)$. The last term considered the case when there is a cut and decompose the initial according to this high-most cut. In that case, we only subtract $1$ in the above part, since the bottom part is just a smaller peeling tree with a distinguished leaf labeled $1$. Notice that when $k=0$, the only possibility to have no leaf labeled $0$ is to have a straight-line from label $1$ to label $r-1$ so that  ${A}_{0}^{(1,r-1)} = 1$ when $r \geq 2$.

\section{``Last Car" decomposition of triangulations}\label{sec:triang}
We now apply our techniques to the case of triangulations. We first describe the encoding of triangulations by their peeling trees which is similar to the case of quadrangulations.
\subsection{Peeling triangulations}\label{sec:peeling-tri}
\begin{figure}[!h]
 \begin{center}
 \includegraphics[width=14cm]{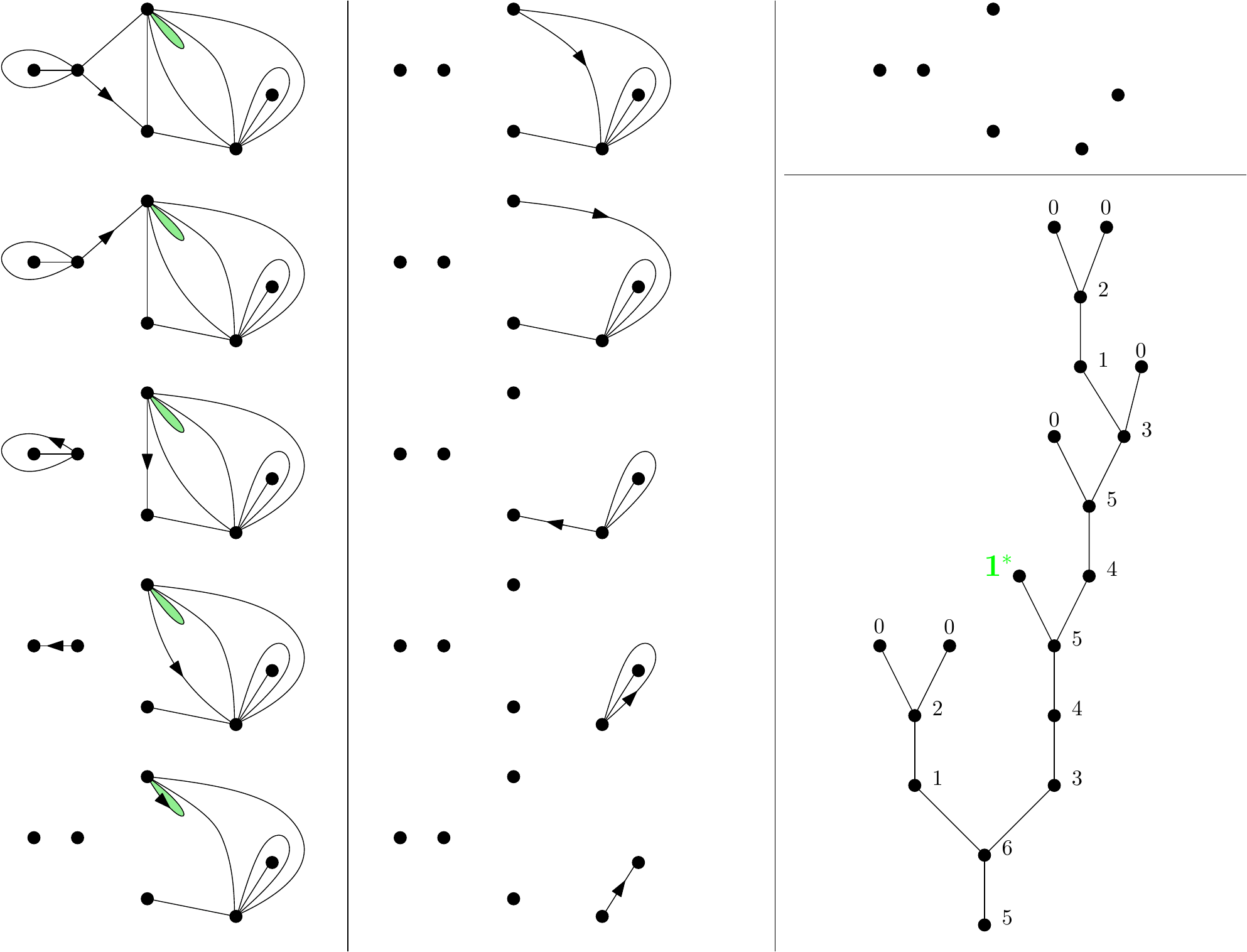}
 \caption{\label{fig:peel-tri} Step-by-step peeling exploration of a triangulation with a boundary of perimeter $p=5$, with $n=6$ vertices and a distinguished loop. This loop in green matches in the tree with the green leaf labeled $1$, which we can put indifferently left or right.}
 \end{center}
 \end{figure}

The peeling technique can be adapted in the case of triangulations. Indeed, we can also remove the edges one-by-one ``à la Tutte" to get recursive equations on the number of triangulations with a boundary. 
As in the case of quadrangulations, two possible events can occur when erasing the root edge of a triangulation with a boundary $p \geq 1$:
\begin{itemize}
\item either the triangulation stays connected, which means that one discovered a new face of the initial triangulation, and one re-roots the triangulation at the left-most edge of the new face. The new triangulation has then a boundary of length $p +1$.
\item or the deletion of the root edge disconnects the triangulation and one gets two triangulations with perimeter $p_1 \geq 0$ and $p_2 \geq 0$ such that $p_1+ p_2 +2 = p$, that one can re-root easily using the two endpoints of the removed edge as in the case of quadrangulations, see Figure \ref{fig:peel-tri}.
\end{itemize}
In that case, we will not consider the half-perimeter of the boundary since the triangulations are not bipartite.
As for quadrangulations, we can encode this edge-by-edge peeling exploration in a peeling tree by recursively labeling the vertices with the perimeter of the boundaries, see Figure \ref{fig:peel-tri}. Note in particular that when the removal of the root edge disconnects the triangulation, we put the component attached to the origin of the root to the left in the peeling tree. To fix ideas, the peeling process gives us the following correspondence between triangulations and trees. \\

\noindent\fbox{\parbox{14cm}{\textbf{Peeling trees of triangulations with one boundary.}\label{encadre-triang}
 To each rooted triangulation with $n$ vertices and with a boundary of perimeter $p$, we can bijectively match its peeling tree
which is a plane tree with labeled vertices with zero, one or two children encoding recursively the perimeter of the boundaries in the peeling exploration. The labeling is such that:
\begin{itemize}
\item the tree has $n$ leaves with label $0$,
\item its root has label $p$, 
\end{itemize} 
\vspace{-0.28cm}
\begin{itemize}[label=$\circ$]
\item  \emph{(local rules)} Each vertex with label $\ell \geq 1$ has either a child with label $\ell+1$ or two children with label $\ell_1$ and $\ell_2$ such that $\ell_1 + \ell_2 +2 = \ell$. See Figure \ref{fig:peel-tri}.
\end{itemize} 
In particular, the local rules imply that it has $3n-p-3$ inner vertices. %
}}
 \\
 
\begin{figure}[!h]
 \begin{center}
 \includegraphics[width=14.5cm]{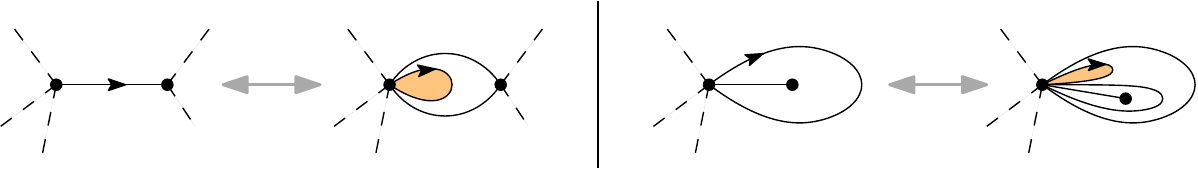}
 \caption{ \label{fig:root-tri} Root transform for triangulations: on the right, the case when the initial root edge is a loop and on the left, the case when the two endpoints of the root edge are different. In both cases, the initial triangulation is on the left and on the right, we obtain a triangulation with a boundary of perimeter $1$.}
 \end{center}
 \end{figure}

As for the case of quadrangulations, it will be convenient to see triangulations of the sphere as triangulations with a boundary. The most natural idea would be to see a triangulation of the sphere as a triangulation with a boundary of perimeter $2$ after unzipping the root edge or $3$ if we see the triangle lying on the right of the root edge as the external face. It implies $[y^2] \mathbf{T}= [y^3] \mathbf{T}$. But since the triangulations are not bipartite, the root edge can be a loop and it will be more convenient to use another root transform. We shall actually view triangulations without boundary as triangulations with a boundary of length $1$ (i.e.\ a loop): we cut along the root edge and “open” it to get a double edge, and then insert a loop inside this double edge at the starting point of the initial root edge to obtain a triangle and root the new triangulation on this loop in clockwise direction so that new triangulation has the $1$-gon to its right, see Figure \ref{fig:root-tri}.  We also obtain $ [y^1] \mathbf{T}=  [y^2] \mathbf{T}= [y^3] \mathbf{T}.$
The choice of the clockwise orientation of the new root edge on the new loop is canonical since we imposed that the maps with a boundary have their distinguished face to their right. 

In fact, we can apply this root-transformation on any distinguished oriented edge even in triangulations which already have a boundary. We will apply our decomposition to peeling trees of triangulations with a boundary and a distinguished loop. When the distinguished loop does not lie on the boundary, we can ``inverse" the root transform and see the distinguished loop as coming from a distinguished oriented edge. To summarize, our decomposition will be based on the following trees. \\
\noindent\fbox{\parbox{14cm}{\textbf{Peeling trees of triangulations with one boundary and a distinguished loop.}
To each rooted triangulation with $n$ vertices and with a boundary of perimeter $p$ and a distinguished loop, we can bijectively match its peeling tree
which is a plane tree with labeled vertices with between $0$ and $2$ children encoding recursively the perimeter of the boundaries in the peeling exploration. The labeling satisfies:
\begin{itemize}
\item the tree has $n$ leaves with label $0$,
\item its root has label $p$, 
\item one inner vertex with label $ \ell \geq 1$ has two children, one of which is a leaf with label $1$ and the other one hase label $\ell - 1$. For this leaf (and only for this leaf !), the planar ordering does not matter, see Figure \ref{fig:peel-tri}.
\end{itemize} 
\vspace{-0.28cm}
\begin{itemize}[label=$\circ$]
\item  \hypertarget{itm:loc-tri}{\emph{(local rules)}} All other vertices with label $\ell \geq 1$ has either a child with label $\ell+1$ or two children with label $\ell_1$ and $\ell_2$ such that $\ell_1 + \ell_2 +2 = \ell$. See Figure \ref{fig:peel-tri}.
\end{itemize} 
In particular, the local rules imply that it has $3n-p-3$ inner vertices. %
}}\\ 

It will then be more convenient to apply our last car decomposition in the branch between this leaf with label $1$ and the root of the tree.

\subsection{Triangulations without boundary}
We first want to enumerate triangulations without boundary and to do this, we apply the above root-transform which gives a bijection between rooted triangulations without boundary and triangulations with a boundary of length $1$ and preserves the number of vertices. 
Instead of distinguishing a vertex to apply our decomposition, we shall this time distinguish an oriented edge, to which we apply another time the above root transform to obtain a triangulation with a boundary of perimeter $1$ (obtained from the transformation of the initial root edge) and with a distinguished loop (obtained from the additional distinguished oriented edge). Note in particular that the oriented distinguished edge can be the root edge, but will be different after applying the root transform for the second time. Therefore the distinguished loop can not lie on the boundary and its edge is part of an internal face so that we can inverse the root transform to recover the initial distinguished oriented edge. 

We build then build its peeling tree to which we can apply our transformation i.e.\ subtract $1$ in the whole branch between the distinguished leaf labeled $1$ and the root of the tree. The transformation works then as in the case of quadrangulations: the local rules (\hyperlink{itm:loc-tri}{$\circ$}) are preserved but two issues may appear.
\begin{itemize}
\item The initial distinguished leaf with label $1$ has a parent with label $k \geq 1$ and a sibling with label $k-1$, so that we just contract  these three vertices into a single marked vertex with label $k-1$.
\item if a vertex labeled $1$ becomes a vertex with label $0$, then we cut the edge just above it; and we remove the initial root which had label $1$.
\end{itemize}

\begin{figure}[!h]
 \begin{center}
 \includegraphics[width=14cm]{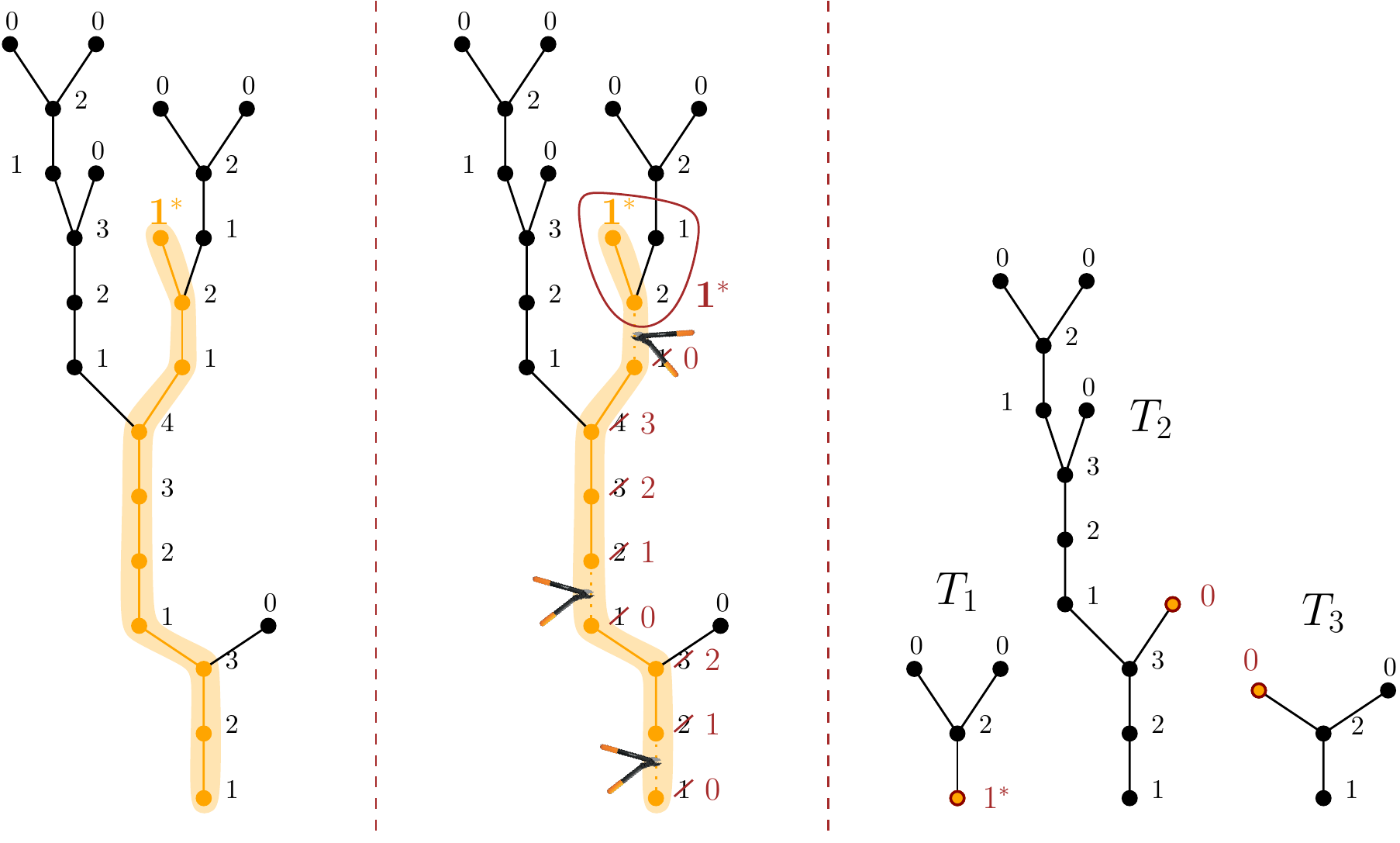}
 \caption{\label{fig:decomp-tri} Example of our decomposition on a tree with $n=6$ leaves labeled $0$ into $k=3$ trees with $n+k-1= 8$ leaves in total. On the left, the initial peeling tree where the distinguished $1$-leaf is displayed in orange as well as the branch between this leaf and the root. After removing $1$ on all labels of the orange branch (in the middle), we need to remove the initial root (or cut just above), cut above the two vertices which get label $0$ and contract the father and sibling of the distinguished leaf into a (marked) vertex labeled $1$. On the right, the resulting sequence of three trees.}
 \end{center}
 \end{figure}

The transformation then takes a tree as described above and maps it into a sequence of $k$ trees with root labeled $1$ and such that the first tree has a distinguished vertex (leaf or not), each possible other tree has a distinguished leaf, and the trees have $n+k-1$ leaves in total, see Figure \ref{fig:decomp-tri}. 

This gives the following differential equation on $ \mathfrak{T}$
\begin{equation*} 
6\mathfrak{T}^{\bullet} - 8 \mathfrak{T} = 4\left(\frac{ \mathfrak{T}^{\bullet}- \mathfrak{T}}{1 - \mathfrak{T}^{ \bullet}/x}\right)
\end{equation*}
which is equivalent to Equation \eqref{eq:tri}. Let us explain  this equality: 
\begin{itemize}
\item On the left, since a triangulation of the $1$-gon has $3n-4$ edges, the function $6\mathfrak{T}^{\bullet} - 8 \mathfrak{T}$ is the generating function of triangulations with a boundary of length $1$ given with an oriented edge.
\item On the right, the first tree of our transformation is a tree which has $n_1$ leaves for some $n_1 \geq 1$ hence $3n_1 - 1 - 3$ inner vertices and $ 4n_1-4$ vertices in total, one of them is distinguished which gives the factor $ 4\mathfrak{T}^{\bullet}- 4\mathfrak{T}$. 
\item The factor $1/(1 - \mathfrak{T}^{ \bullet}/x)$ gives the potential other trees which are marked at a leaf. 
\end{itemize}
 By identifying the coefficients in the above equation, a straightforward computation gives Equation \eqref{eq:tri-rec}, which characterizes entry A002005 in Sloane online encyclopedia for integer sequences.

\subsection{Triangulations with one boundary}
As in the case of quadrangulations, the transformation which we described above can be adapted in the case of triangulations with a boundary of perimeter $p \geq 2$. More precisely, we consider peeling trees of triangulations with a boundary of length $p \geq 2$ and one distinguished loop that we described in Section \ref{sec:peeling-tri}. Such peeling trees have a distinguished leaf with label $1$ coming from  the distinguished loop for which we recall its local rule: if its father has label $\ell$, its sibling has label $\ell-1$. We can apply now our ``last car"-transformation on the branch between this leaf labeled $1$ and the root of the peeling tree. The only difference with that of peeling trees starting from $1$ is that we do not remove the root of the trees starting from $p \geq 2$. Our last car decomposition of a peeling tree starting from $p$ with $n$ leaves labeled $0$ produces then a sequence of $k$ peeling trees of triangulations for some $k \geq 0$ such that: 
\begin{itemize}
\item the first tree is marked at a vertex (leaf or not)
\item all eventual following trees have a marked leaf
\item all but the last tree have a root label $1$ and the last tree has a root labeled $p - 1 \geq 1$
\item the trees have $n+k -1$ leaves labeled $0$ in total.
\end{itemize}
Conversely, given the sequence, we can recover the initial tree by stacking and gluing the trees and add $1$ in the appropriate branch so that it is really a bijection, which gives Equation \eqref{eq:tri1bord}, which we recall here:
$$6 \mathbf{T}^{\bullet} - 2 y \partial_{y}\mathbf{T} - 6\mathbf{T} + y \partial_{y} (y \mathbf{T}) =  \frac{4y}{x} \mathbf{T}^{\bullet} \left(\frac{ \mathfrak{T}^{\bullet}- \mathfrak{T}}{1 - \mathfrak{T}^{\bullet}/x}\right) + y \left(4\mathbf{T}^{\bullet}- 3\mathbf{T}- y\partial_{y}\mathbf{T}\right). $$
where $ \mathbf{T}^{\bullet} = x \partial_{x}\mathbf{T}$. 
Let us explain this equation in more details: 
\begin{itemize}
\item The left-hand side enumerates triangulations with a boundary of length $p \geq 2$ and one distinguished loop with a weight $x$ per vertex and $y$ per boundary length. If the loop is on the boundary, then we can map the triangulation to a triangulation with perimeter $p -1$ and a distinguished vertex on the boundary where the initial loop is attached. Therefore, this type of triangulation are enumerated by the factor $y \partial_{y} (y \mathbf{T})$. When the loop is not on the boundary, then this edge is really an edge of an inner triangle so that we can perform the converse of the root transform and obtain a distinguished oriented edge. Since there are $3n-p-3$ (unoriented) edges in a triangulation with $n$ vertices and boundary of length $p$, these are enumerated by the term $6 \mathbf{T}^{\bullet} - 2 y \partial_{y}\mathbf{T} - 6\mathbf{T}$, which explained the left-hand side term. 
\item On the right, the first term corresponds to the case when our transformation gives more than a tree (and in particular, the first and the last tree are not the same), including the case of triangulation with boundary $1$  where the last tree after transformation is the tree with one vertex labeled $0$. The last tree is in that case a tree with root label $p-1$ and a marked leaf, which corresponds to the factor $y \mathbf{T}^{ \bullet}$. The first tree has a distinguished vertex and its root has label $1$, hence the factor $4 \mathfrak{T}^{\bullet} - 4 \mathfrak{T}$, and the division by $(1 - \mathfrak{T}^{\bullet}/x)$ stands for the other possible cuts. There is no factor $2$ in that case since the discovery of the marked loop can be indifferently put left or right in the tree, but a factor $1/x$ to obtain the appropriate total number of leaves.
\item The last term corresponds to the triangulations where we get only one tree when applying our transformation on its peeling tree. In that case the resulting tree has $n-1$ leaves, a root labeled $p-1$ hence $3(n-1)-(p - 1)-3$ inner vertices in the tree and $ 4(n-1) - 3 - (p-1) $ vertices in total, one of them is distinguished explaining the factor $4\mathbf{T}^{\bullet}- 3\mathbf{T}- y\partial_{y}\mathbf{T}$.
\end{itemize}

\subsection{Triangulations with two boundaries}

As in the case of quadrangulations, we are also able to enumerate triangulations with two boundaries. Considering the second boundary as a distinguished face, we can do the same decomposition according to the discovery of this different face. The discovery of this distinguished face corresponds in the peeling tree to a transition from a vertex labeled $r \geq 1$ to a vertex labeled $r+q-2$ since we consider here the perimeter (and not the half-perimeter). \\
We introduce for this purpose ${B}_{n}^{(p,r)}$ the number of peeling trees of triangulations with $n$ leaves $0$ and a distinguished leaf with label $r$ and boundary of length $p$. When $r=0$, we set $ B_n^{(p,0)}= (n+1)T_n^{(p)}$ when $n \geq 1$. The last car decomposition leads to the following recursive equations:
\begin{eqnarray*}
\forall r \geq 1, \forall p \geq 1, \forall n \geq 1, &{B}_{n}^{(p,r)} &= {B}_{n}^{( \max(p-1,1) , r-1)} + \sum_{k=0}^{n-1} {B}_{k}^{(1,r-1)} {B}_{n-k}^{(p,1)}.
\end{eqnarray*}
This equation are very similar to  Equation  \eqref{eq:recTreeQuadR}. We only need to adapt the initial conditions in the $B_{n}^{(p,0)}$'s. 
Given the $ {B}_{n}^{(p,r)}$'s,  we are now able to deduce the number $T_n^{(p,q)}$ of rooted triangulations with $n$ vertices and with a boundary of length $p$ and a distinguished face of length $q$ from this equation
$$ {T}_{n}^{(p,q)} = \sum_{r=1}^{2n-q}\sum_{k = 0}^{n- \lfloor \frac{r+q}{2}\rfloor+1} {B}_{k}^{(p,r)} T_{n-k}^{(r+q-2)},$$
which is the analog of  Equation \eqref{eq:rec2bords} in the case of quadrangulations with two boundaries. The only changes are that we replaced the initial conditions and adjusted the bounds of the sums. Notice that we already enumerated the the triangulation with one boundary and a distinguished loop (the $T_n^{(p,1)}$'s) by $ \mathbf{T}^{\bullet} - 2 y \partial_{y}\mathbf{T} - 6\mathbf{T} + y \partial_{y} (y \mathbf{T}) $ the left-hand side term of~\eqref{eq:tri1bord}.

\section{Comments and perspectives}
We mention here a few possible developments of this work.
 \paragraph{Other models of maps.} We applied here a ``last car decomposition" in the case of quadrangulations and triangulations with zero, one or two boundaries. We believe that this decomposition can be adapted for other models of planar maps such as $p$-angulations (at least for $p$ even) or (bipartite) maps with Boltzmann weights.

\paragraph{Solving equations.}
We gave here new equations which characterize the enumeration of a certain type of maps. Some of them were already explicitly enumerated (quadrangulations with zero, one or two boundaries, triangulations without boundary...) so it can be easily checked that their generating functions satisfy our equations. But we may wonder if we can recover those coefficients directly by solving our equations explicitly. This is also a relevant question in the cases when no explicit formula is known: can we extract from these equations explicit formulas for the coefficients? 

\paragraph{Decomposition of maps.}
Here we gave a decomposition of the peeling tree into smaller peeling trees i.e.\ peeling trees of ``smaller" maps.  We do not know if this decomposition can be easily interpreted on the maps. In particular, we choose a specific way to reroot the maps in the peeling exploration, but there are many ways to do it by choosing different \emph{peeling algorithm}. There may be a choice of peeling algorithm for which the transformation on the maps is ``natural". We have no hope that it is local but it may have similarities with the cut and slice operation of Louf \cite{louf2019new}.

\bibliographystyle{siam}
\bibliography{/Users/contat/Dropbox/Articles/biblio}

\begin{thebibliography}{10}

\bibitem{BMJ06}
{\sc M.~Bousquet-M{\'e}lou and A.~Jehanne}, {\em Polynomial equations with one
  catalytic variable, algebraic series and map enumeration}, J. Combin. Theory
  Ser. B, 96 (2006), pp.~623--672.

\bibitem{BIPZ78}
{\sc E.~Br{\'e}zin, C.~Itykson, G.~Parisi, and J.-B. Zuber}, {\em Planar
  diagrams}, Comm. Math. Phys.,  (1978).

\bibitem{chen2021enumeration}
{\sc L.~Chen}, {\em Enumeration of fully parked trees}, arXiv preprint
  arXiv:2103.15770,  (2021).

\bibitem{chen2021parking}
{\sc Q.~Chen and C.~Goldschmidt}, {\em Parking on a random rooted plane tree},
  Bernoulli, 27 (2021), pp.~93--106.

\bibitem{contat2020sharpness}
{\sc A.~Contat}, {\em Sharpness of the phase transition for parking on random
  trees}, Random Structures \& Algorithms,  (2020).

\bibitem{contat2021parking}
{\sc A.~Contat and N.~Curien}, {\em Parking on {C}ayley trees \& {F}rozen
  {E}rd{\H{o}}s-{R}{\'e}nyi}, arXiv preprint arXiv:2107.02116,  (2021).

\bibitem{CoriSchaefferDescription}
{\sc R.~Cori and G.~Schaeffer}, {\em Description trees and {T}utte formulas},
  Theoretical Computer Science, 292 (2003), pp.~165--183.

\bibitem{CV81}
{\sc R.~Cori and B.~Vauquelin}, {\em Planar maps are well labeled trees},
  Canad. J. Math., 33 (1981), pp.~1023--1042.

\bibitem{CurStFlour}
{\sc N.~Curien}, {\em Peeling random planar maps, Saint-Flour course 2019},
  https://www.imo.universite-paris-saclay.fr/$\sim$curien/.

\bibitem{CH19}
{\sc N.~Curien and O.~H\'enard}, {\em The phase transition for parking on
  {G}alton-{W}atson trees}, arXiv:1912.06012,  (2019).

\bibitem{fang2021bijective}
{\sc W.~Fang}, {\em Bijective link between {C}hapoton's new intervals and
  bipartite planar maps}, European Journal of Combinatorics, 97 (2021),
  p.~103382.

\bibitem{GP19}
{\sc C.~Goldschmidt and M.~Przykucki}, {\em Parking on a random tree},
  Combinatorics, Probability and Computing, 28 (2019), pp.~23--45.

\bibitem{GJ08}
{\sc I.~P. Goulden and D.~M. Jackson}, {\em The {KP} hierarchy, branched
  covers, and triangulations}, Adv. Math., 219 (2008), pp.~932--951.

\bibitem{konheim1966occupancy}
{\sc A.~G. Konheim and B.~Weiss}, {\em An occupancy discipline and
  applications}, SIAM Journal on Applied Mathematics, 14 (1966),
  pp.~1266--1274.

\bibitem{LaP16}
{\sc M.-L. Lackner and A.~Panholzer}, {\em Parking functions for mappings},
  Journal of Combinatorial Theory, Series A, 142 (2016), pp.~1 -- 28.

\bibitem{louf2019new}
{\sc B.~Louf}, {\em A new family of bijections for planar maps}, Journal of
  Combinatorial Theory, Series A, 168 (2019), pp.~374--395.

\bibitem{panholzer2020combinatorial}
{\sc A.~Panholzer}, {\em A combinatorial approach for discrete car parking on
  random labelled trees}, Journal of Combinatorial Theory, Series A, 173
  (2020), p.~105233.

\bibitem{schaeffer1997bijective}
{\sc G.~Schaeffer}, {\em Bijective census and random generation of {E}ulerian
  planar maps with prescribed vertex degrees}, Electronic Journal of
  Combinatorics, 4 (1997), p.~20.

\bibitem{tHo74}
{\sc G.~'t~Hooft}, {\em A planar diagram theory for strong interactions.},
  Nuclear Physics B, 72 (1974), pp.~461--473.

\bibitem{Tut62c}
{\sc W.~T. Tutte}, {\em A census of {H}amiltonian polygons}, Canad. J. Math.,
  14 (1962), pp.~402--417.

\bibitem{Tut62}
\leavevmode\vrule height 2pt depth -1.6pt width 23pt, {\em A census of planar
  triangulations}, Canad. J. Math., 14 (1962), pp.~21--38.

\bibitem{Tut62b}
\leavevmode\vrule height 2pt depth -1.6pt width 23pt, {\em A census of
  slicings}, Canad. J. Math., 14 (1962), pp.~708--722.

\bibitem{Tut63}
\leavevmode\vrule height 2pt depth -1.6pt width 23pt, {\em A census of planar
  maps}, Canad. J. Math., 15 (1963), pp.~249--271.

\end{thebibliography}
\end{document}